\documentclass[11pt,reqno]{article}
\usepackage{amssymb,eucal}
\usepackage{amsmath}
\usepackage{color}

\newtheorem{thm}{Theorem}[section]

\newtheorem{lem}{Lemma}[section]
\newtheorem{rem}{Remark}[section]

\newtheorem{expl}{Example}[section]

\makeatletter
   
          \@addtoreset{equation}{section}

\begin{document}

\title{
Schnol's theorem and spectral properties of massless  Dirac
operators with scalar potentials.
}

\author{Karl Michael Schmidt$^1$~~and Tomio Umeda$^2\!$
\footnote{Supported by
 the Japan Society for the Promotion of Science
     ``Grant-in-Aid for Scientific Research'' (C)
    No.  21540193.}}

\date{}
\maketitle

\begin{quote}
\begin{itemize}
\item[$^1$] School of Mathematics,  Cardiff University, Senghennydd Road,\\
Cardiff CF24 4AG, Wales, UK.\\
E-mail: {\tt schmidtkm@cardiff.ac.uk}
\item[$^2$] Department of Mathematical Sciences, University of Hyogo, Shosha,\\ Himeji 671-2201, Japan \\
E-mail: {\tt umeda@sci.u-hyogo.ac.jp}
\end{itemize}
\end{quote}

\begin{abstract}
\noindent
The spectra of massless Dirac operators are of essential interest
e.g. for the electronic properties of graphene, but fundamental questions such
as the existence of spectral gaps remain open.
We show that the eigenvalues of massless Dirac operators with suitable
real-valued potentials lie inside small sets easily characterised in terms of
properties of the potentials, and we prove a Schnol'-type theorem relating
spectral points to polynomial boundedness of solutions of the Dirac equation.
Moreover, we show that, under minimal hypotheses which leave the potential
essentially unrestrained in large parts of space, the spectrum of the massless
Dirac operator covers the whole real line;
in particular, this will be the case if the potential is nearly constant
in a sequence of regions.
\end{abstract}


\smallskip


\vspace{16pt}

\section{Introduction}

\noindent

The Dirac operators we shall consider  in this paper are 
\begin{equation}  \label{eq:do2d}
H_2= -i \,\sigma \cdot \!\nabla + q(x)  \;\; \mbox{ in } \;
{\mathsf L}^2({\mathbb R}^2; {\mathbb C}^2)
\end{equation}
and 
\begin{equation}   \label{eq:do3d}
H_3= -i \, \alpha  \cdot \! \nabla + q(x)  \;\; \mbox{ in } 
\;
{\mathsf L}^2({\mathbb R}^3; {\mathbb C}^4),
\end{equation}
where  $\sigma=(\sigma_1, \, \sigma_2)$ 
and 
$\alpha=(\alpha_1, \, \alpha_2, \, \alpha_3)$
are given as follows:
\begin{equation*}\label{eqn:1-3}
\sigma_1 =
\begin{pmatrix}
0&1 \\ 1& 0
\end{pmatrix}, \,\,\,
\sigma_2 =
\begin{pmatrix}
0& -i  \\ i&0
\end{pmatrix}
\end{equation*}
and
\begin{equation*}\label{eqn:1-4}
\alpha_j =
\begin{pmatrix}
0&\sigma_j \\ \sigma_j& 0
\end{pmatrix}
\quad (\, j\in \{1, \, 2, \, 3\})
 \,\,\,
\mbox{ with } \,\,\,
\sigma_3 =
\begin{pmatrix}
1&0 \\ 0&-1
\end{pmatrix}\!.
\end{equation*}

\vspace{4pt}
\noindent
The dot products are to be read as
\begin{equation*}
\sigma \cdot \!\nabla= \sigma_1 \frac{\partial}{\partial x_1} 
+\sigma_2  \frac{\partial}{\partial x_2} 
\end{equation*}
in (\ref{eq:do2d})
and 
\begin{equation*}
\alpha \cdot \!\nabla= \alpha_1 \frac{\partial}{\partial x_1} 
+\alpha_2  \frac{\partial}{\partial x_2} 
+  \alpha_3  \frac{\partial}{\partial x_3} 
\end{equation*}
in (\ref{eq:do3d}).
The potential 
$q$ is a real-valued function on ${\mathbb R}^d$, 
where $d=2$ or $d=3$, respectively.
The operators $H_2$ and $H_3$ differ from the standard Dirac operator
in that they lack
a mass term, usually represented by an additional anti-commuting
matrix: $\sigma_3$ for the two-dimensional case and 
\begin{equation*}
\beta= 
\begin{pmatrix}
I & 0  \\  0 & -I 
\end{pmatrix}
\end{equation*}
for the three-dimensional case, where $I$ is a $2\times2$ identity matrix.

The purpose of the present paper
 is twofold.
Firstly, we establish Schnol's theorem for $H_d$, $d=2$ and $3$,
under minimal assumptions on $q$.   
Schnol's theorem
for Schr\"odinger operators is well-known (cf. 
 \cite[p.21]{CFKS}); it asserts that
an energy with polynomially bounded
eigensolution  belongs to the spectrum
of the Schr\"odinger operator. 
In  this context, 
we would like to mention some recent works 
on Schnol's theorem for 
generators of Dirichlet forms,
cf. \cite{BMLS}, \cite{HK} and \cite{LSV}.
To our knowledge, however, the present paper is the first to
establish  Schnol's theorem for Dirac operators.
Secondly, we shall
 show that
 $\sigma(H_d) = \mathbb R$
under minimal assumptions on $q$ as before.
We shall not require any restriction on the growth 
or decay of
 the potential $q$ at infinity.
\vspace{10pt}

The study of the spectrum of massless Dirac operators in the two- and
three-dimensional case is intriguing, as the behaviour of these operators
differs fundamentally from the more familiar cases of Schr\"odinger
operators and Dirac operators with mass. The 
spectrum of the one-dimensional massless 
Dirac operator 
\begin{equation} \label{eq:do1d}
H_1 =- i \sigma_2 \frac{d}{dx} + q(x)
\;\; \mbox{ in } \;
{\mathsf L}^2({\mathbb R}; {\mathbb C}^2)
\end{equation}
covers the whole real axis and  is purely absolutely continuous
whenever $q \in {\mathsf L}^1_{loc}(\mathbb R, \, \mathbb R)$.
This surprising fact  was first pointed out 
 by one of the authors in  \cite{KMS-1}.
By separation in spherical polar coordinates, this
result also implies that 
$\sigma(H_d)= \mathbb R$
if $q$ is rotationally symmetric.
However, the situation is by no means clear in the more general
higher-dimensional case. The two-dimensional massless Dirac
operator is of particular interest because it governs electron
transport in graphene, so its spectral properties will have a
direct impact on the conductivity and potential use in electronic
applications.
It is known that total reflection of the quantum wave at a
straight-edged potential step may occur \cite{BTB}, and initially
there was some hope to capture bound states in localised quantum
dots (see \cite{BTB}, Fig. 1(b)).
However, this is impossible due to a result of \cite{KOY} which,
in particular, implies that a compactly supported potential cannot
generate eigenvalues (see \cite{KOY}, Ex.~6.1).
Furthermore,
it is believed that
the energy spectrum of graphene, 
irrespective of potential applied,
has no bandgap (zero bandgap); see \cite{BTB}, \cite{DSDLA},
\cite{NG}.
This question remains open, but from the results mentioned above it is
clear that spectral phenomena such as gaps or eigenvalues will, if
they occur at all, require potentials of a fairly complex global structure.
Recently, the properties of disordered graphene have attracted much attention,
and it is known
\cite{Peres} that the sources of disorder vary and can be 
described by various types of potentials.
The dominant source of disorder is still under debate according to
\cite{Zhaoetal}. 
Under these circumstances,
 it is natural, 
from the mathematical point of view, 
to investigate spectral properties of $H_d$ and, in particular,
to make an attempt to show that
 $\sigma(H_d) = \mathbb R$, under minimal assumptions on the potential $q$.

An announcement of the present paper can be found in
\cite{SU}.

\vspace{10pt}

\section{Embedded eigenvalues and the absolutely continuous spectrum}
In contrast to the case of the one-dimensional Dirac operator $H_1$,
the spectra of 
$H_2$ and $H_3$ are 
not
always purely absolutely
continuous 
regardless of $q$. 
Actually, in the three-dimensional case,
 we have an example 
of $q$ which gives rise to a zero mode of $H_3$, 
\textit{i.e.} an example of
$q$ for which $H_3$ has the embedded eigenvalue $0$.

\begin{expl}  \label{ex:3}
Let
$q(x)= -3/ \langle x \rangle^2$, where
$\langle x \rangle= \sqrt{1 + |x|^2}$.
Then  there exists
a unique self-adjoint realization of $H_3$ in
${\mathsf L}^2({\mathbb R}^3; {\mathbb C}^4)$
with 
$\mbox{\rm Dom}(H_3)={\mathsf H}^1({\mathbb R}^3; {\mathbb C}^4)$, 
the Sobolev space
of order $1$. 
If one puts
\begin{equation*}  \label{eq:zeromd3}
f(x)= \langle x \rangle^{-3} (I_4 + i \alpha \cdot x) \phi_0
\end{equation*}
 with 
$\phi_0$ a unit vector in ${\mathbb C}^4$, 
then a direct calculation shows that $H_3f=0$.
Since $f \in {\mathsf H}^1({\mathbb R}^3; {\mathbb C}^4)$, 
this implies that $0 \in \sigma_{\rm p}(H_3)$.
Thus $H_3$ has a zero mode. 
As $\lim\limits_{|x|\rightarrow\infty} q(x) = 0$, a simple singular sequence
argument shows that $\sigma(H_3)= \mathbb R$. Hence the energy $0$ is 
an embedded eigenvalue of $H_3$.
\end{expl}

We would like to mention that
the potential $q$ and the zero mode $f$  in Example \ref{ex:3}
were motivated by \cite{LY}.

The analogous two-dimensional
construction in Example \ref{ex:2} below gives
a zero resonance of $H_2$,  not a zero mode of $H_2$.
We do not know if the potential $q$ 
in  Example \ref{ex:2} gives rise to a zero mode of  $H_2$.
However, zero modes of $H_2$ are known to occur with compactly supported
rotationally symmetric potentials, see Theorem 3 of \cite{KMS-2}.

\begin{expl}  \label{ex:2}
Let
$q(x)= -2/ \langle x \rangle^2$.
Then there exists
a unique self-adjoint realization of $H_2$ in
${\mathsf L}^2({\mathbb R}^2; {\mathbb C}^2)$
with 
$\mbox{\rm Dom}(H_2)
=
{\mathsf H}^1({\mathbb R}^2; {\mathbb C}^2)$,
and $\sigma(H_2)= \mathbb R$. 
If
\begin{equation*}
\psi(x) =\langle x \rangle^{-2} (I_2 + i \sigma \cdot x) \phi_0,
\end{equation*} 
$\phi_0$ a unit vector in ${\mathbb C}^2$,
then one sees that $H_2 \psi =0$.
However, it is clear that
 $\psi \not \in {\mathsf L}^2({\mathbb R}^2; {\mathbb C}^2)$.
Therefore, $\psi$ is not a zero mode of $H_2$. 
On the other hand, one finds that
$\psi \in 
{\mathsf L}^{2}_{-s}({\mathbb R}^2; {\mathbb C}^2)$ 
for any $s>0$, 
where
\begin{equation*}
{\mathsf L}^{2}_{-s}({\mathbb R}^2; {\mathbb C}^2)
=
\big\{ \, \varphi \;  \big| \;  
 \Vert \langle x \rangle^{-s}\varphi \Vert_{L^2} < \infty \big\}.
\end{equation*}
This means that $\psi$ is a zero resonance of $H_2$.
\end{expl}

It is not 
an easy task  to clarify whether $H_d$ has
embedded eigenvalues for general potentials. 
However,
we have a good control of the embedded 
eigenvalues of $H_d$ if
$q(x)$ is rotationally symmetric. 
To formulate a result, 
we need to introduce the definition of 
the limit range  $\mathcal R_{\infty}(q)$ of $q$:
\begin{equation*}
\mathcal R_{\infty}(q) = \bigcap_{r>0} 
\overline{\,\big\{ 
\, q(x)  \, \big| \, |x|\ge r \big\} },
\end{equation*}
where $\overline{A}$ denotes the closure of a subset of 
$A \subset \mathbb R$.

\begin{thm}[\textbf{Schmidt\cite{KMS-2}}]  \label{th:KMS}
Let $q(x)=\eta(|x|)$ and let $\eta \in 
{\mathsf L}^1_{loc}(0, \, \infty)$.
Suppose that there exists a real number
$\lambda \in \mathbb R \setminus \mathcal R_{\infty}(q)$
such that 
\begin{equation} \label{eq:BV}
\frac{1}{\, r(\lambda - \eta (r) )-1 \,} \in BV(r_0, \infty) 
\end{equation}
for some $r_0 >0$,
where $BV(r_0, \infty)$ denotes
the set of functions of bounded variations on 
the interval $(r_0, \infty)$.
Then $\sigma_p(H_d) \subset \mathcal R_{\infty}(q)$.
\end{thm}

\vspace{3pt}

Theorem \ref{th:KMS} is a direct consequence of
\cite[Corollary 1]{KMS-2}.
One should note that under the assumption that 
$\eta \in {\mathsf L}^1_{loc}(0, \, \infty)$
there exists a distinguished  self-adjoint
realization of $H_d$, see Propositions 2 and 3 of \cite{KMS-2}.

Note that the condition (\ref{eq:BV}) just fails in the radially periodic case,
i.e. if $\eta(r+p) = \eta(r)$ $(r \ge 0)$ with period $p > 0$.
Indeed, when we take $r, s \in [0, \, p]$ and $n \in {\mathbb N}$, then for $\lambda$ not
in the (limit) range of $\eta$,
\begin{align*}
&\frac{1}{(r+np)(\lambda - \eta(r+np)) - 1} - \frac{1}{(s+np)(\lambda - \eta(s+np)) - 1}
\\
&\qquad\sim \frac{1}{(r+np)(\lambda - \eta(r))} - \frac{1}{(s+np)(\lambda - \eta(s))}
\\
&\qquad= \frac{sr}{(s+np)(r+np)} \left(\frac{1}{r\,(\lambda - \eta(r))} - \frac{1}{s\,(\lambda - \eta(s))}\right)
\\
&\qquad\qquad+ \frac{np}{(s+np)(r+np)}\left(\frac{1}{\lambda - \eta(r)} - \frac{1}{\lambda - \eta(s)}\right)
\end{align*}
as $n \rightarrow \infty$, so the total variation in the $n$th period interval is
\begin{equation*}
\mathop{\textrm{Var}}\limits_{r \in [np, \,(n+1)p]} \frac{1}{r(\lambda - \eta(r)) - 1)}
\sim \frac{1}{np}\,\mathop{\textrm{Var}}\limits_{[0, \, p]} \frac{1}{\lambda - \eta}
\qquad
(n\rightarrow\infty).
\end{equation*}
In fact, the limit range of the potential plays no role at all in the
radially periodic case, as our following result shows.

\begin{thm} \label{th:radper}
Let $q(x) = \eta(|x|)$ and let $\eta \in L^1_{loc}(0,\infty)$ be $p$-periodic.
Let $\hat\eta := \frac{1}{p} \int_0^p \eta$.
Then $H_d$ has purely absolutely continuous spectrum in ${\mathbb R}\setminus (\frac{\pi}{p}{\mathbb Z}
+ \hat\eta)$.
\end{thm}

\medskip
\noindent
{\bf Proof.}
By a suitable shift of the spectral parameter, we may assume without loss of generality
that $\hat\eta = 0$.

By separation of variables in polar coordinates (see e.g. \cite{Weidmann}, Appendix to
Section 1),
\begin{equation*}
H_d \cong \bigoplus_{k\in J_d} -i\sigma_2 \frac{d}{dr} + \eta(r) + \sigma_1 \frac{k}{r},
\end{equation*}
where the index set $J_d = {\mathbb Z}\setminus\{0\}$ if $d = 3$ and $J_d = {\mathbb Z} + \frac{1}{2}$
if $d = 2$.
Hence it is sufficient to show that each of the half-line Dirac operators with the
angular momentum term $\sigma_1 \frac{k}{r}$ has purely absolutely continuous spectrum
in ${\mathbb R}\setminus \frac{\pi}{p}{\mathbb Z}$.

It follows from Gilbert-Pearson subordinacy theory (\cite{GP}, \cite{Gilbert};
see also \cite{KMS-3}) that it is sufficient
for this purpose to show that all solutions of the eigenvalue equation
\begin{equation}\label{eq:raddir}
-i\sigma_2 \frac{d}{dr}\,u(r) + \eta(r)\,u(r) + \sigma_1 \frac{k}{r}\,u(r) = \lambda\,u(r)
\end{equation}
are bounded at infinity if $\lambda \notin \frac{\pi}{p}{\mathbb Z}$.

Let $\varepsilon > 0$; we shall prove the boundedness of all solutions of (\ref{eq:raddir})
for large $r$ and for all
$\lambda \in \Lambda := {\mathbb R} \setminus \bigcup_{n \in {\mathbb Z}} (\frac{n\pi}{p} - \varepsilon,
\frac{n\pi}{p} + \varepsilon)$ by adapting an idea of Stolz \cite{Stolz}; see also
\cite[Theorem 5.2.1]{BES}.

Let $Q(r) := \int_0^r \eta$ $(r \ge 0)$; then $Q$ is $p$-periodic and $Q(0) = 0$.
For $j \in {\mathbb N}$, the matrix-valued function 
$\Phi_j(r) := \exp\big[-i\sigma_2 \big\{ Q(r) - \lambda (r-(j-1)p)\big\}\big]$ $(r \ge 0)$ satisfies
the unperturbed, periodic differential equation
\begin{equation*}
-i\sigma_2 \frac{d}{dr}\,\Phi_j(r) + \eta(r)\,\Phi_j(r) = \lambda\,\Phi_j(r),
\end{equation*}
and $\Phi_j((j-1)p) = I$.
$M(\lambda) := \Phi_j(jp) = I \cos \lambda p + i \sigma_2 \sin \lambda p$ is the monodromy
matrix  (cf. \cite[p.5]{BES}) of the periodic equation and $D(\lambda) := \mathop{\textrm{tr}} M(\lambda) = 2 \cos \lambda p$
its discriminant (cf. \cite[p.9]{BES}).
Clearly there exists $\delta > 0$ such that $|D(\lambda)| \le 2 - 2\delta$ $(\lambda \in
\Lambda)$.

By the variation of constants method (cf. \cite[p.3]{BES}), 
we can find an integral equation for the matrix-valued
solution $\Psi_j$ of (\ref{eq:raddir}) such that $\Psi_j((j-1)p) = I$,
\begin{equation} \label{eq:inteq}
\Psi_j(r) = \Phi_j(r) - \Phi_j(r) \int_{(j-1)p}^r \Phi_j(s)^{-1} \sigma_3 \frac{k}{s}\,
  \Psi_j(s)\,ds
\qquad
(r \ge (j-1)p).
\end{equation}
Using the fact that $\Phi_j$ is always unitary, we hence obtain the estimate for the
matrix operator norm
\begin{equation*}
|\Psi_j(r)| \le 1 + \int_{(j-1)p}^r \frac{|k|}{s}\,|\Psi_j(s)|\,ds
\qquad
(r \ge (j-1)p).
\end{equation*}
By Gronwall's inequality, it follows that
\begin{equation*}
|\Psi_j(r)| \le \exp\left(|k| \log\frac{r}{(j-1)p}\right)
= \frac{r^{|k|}}{(j-1)^{|k|} p^{|k|}},
\end{equation*}
and hence by (\ref{eq:inteq}), for $(j-1)p \le r \le jp$,
\begin{align}
|\Psi_j(r) - \Phi_j(r)| &\le \int_{(j-1)p}^r \frac{|k|}{s}\,|\Psi_j(s)|\,ds
\le \frac{r^{|k|} - (j-1)^{|k|} p^{|k|}}{(j-1)^{|k|} p^{|k|}} \nonumber \\
&\le \frac{j^{|k|} - (j-1)^{|k|}}{(j-1)^{|k|}}
= \big( 1 + \frac{1}{j-1} \big)^{|k|} - 1
\rightarrow 0 \qquad (j\rightarrow\infty). \label{eq:psiphi}
\end{align}
In particular, the matrices $M_j(\lambda) := \Psi_j(jp)$ satisfy
\begin{equation}\label{eq:monocon}
\lim\limits_{j\rightarrow\infty} |M_j(\lambda) - M(\lambda)| = 0
\end{equation}
uniformly in
$\lambda \in {\mathbb R}$.
This implies that $D_j(\lambda) := \mathop{\textrm{tr}} M_j(\lambda) \rightarrow D(\lambda)$ uniformly
in $\lambda$ as $j\rightarrow\infty$.
Thus there is $J \in {\mathbb N}$ such that $|D_j(\lambda)| \le 2 - \delta$ for all $j \ge J$
and $\lambda \in \Lambda$.

For such $j$ and $\lambda$, the matrices $M_j(\lambda)$ have complex conjugate eigenvalues
$\mu_j(\lambda), \overline{\mu_j(\lambda)}$, $|\mu_j(\lambda)| = 1$.
(See \cite[{\it Case 3} \, in p. 10]{BES} in conjunction with the fact that
 $\text{det}\Psi_j(r)=1$, which is obtained from
 \cite[Liouville's formula in p.3]{BES}).
A suitable matrix of eigenvectors can be written as
\begin{equation*}
E_j(\lambda) := \left(\begin{matrix}
 \mu_j(\lambda) - (M_j(\lambda))_{22} & \overline{\mu_j(\lambda)} - (M_j(\lambda))_{22}\\
 (M_j(\lambda))_{21} & (M_j(\lambda))_{21}
  \end{matrix}\right);
\end{equation*}
in the limit $j \rightarrow \infty$, this converges uniformly in $\lambda \in \Lambda$
to a corresponding matrix
$E(\lambda)$ of eigenvectors of $M(\lambda)$ in view of (\ref{eq:monocon}).

Now consider the matrix-valued solution $\Psi_J$. For $n \ge J$ (omitting the variable
$\lambda$ for brevity),
\begin{align*}
\Psi_J(np) &= M_n M_{n-1} \cdots M_J \\
 &= E_n {\textstyle\left(\begin{matrix}\mu_n & 0 \\ 0 & \overline{\mu_n}\end{matrix}\right)} E_n^{-1}
 E_{n-1} {\textstyle \left(\begin{matrix}\mu_{n-1} & 0 \\ 0 & \overline{\mu_{n-1}}\end{matrix}\right)} E_{n-1}^{-1} \cdots
 E_J {\textstyle \left(\begin{matrix}\mu_J & 0 \\ 0 & \overline{\mu_J}\end{matrix}\right)} E_J^{-1}.
\end{align*}
Hence the matrix operator norm can be estimated as
\begin{equation*}
|\Psi_J(np)| \le |E_n|\,|E_n^{-1} E_{n-1}|\,|E_{n-1}^{-1} E_{n-2}| \cdots |E_{J+1}^{-1} E_J|\, |E_J^{-1}|.
\end{equation*}

In order to estimate $|E_j^{-1} E_{j-1}|$, we again solve (\ref{eq:raddir}) on the
interval $[(j-1)p, jp]$ by variation of constants, but this time using $\Psi_{j-1}(r-p)$
as a reference instead of $\Phi_j(r)$; this gives
\begin{equation*}
\Psi_j(r) = \Psi_{j-1}(r-p) + \Psi_{j-1}(r-p) \int_{(j-1)p}^r \Psi_{j-1}(s-p)^{-1}
  \sigma_3\, \frac{kp}{s(s-p)}\,\Psi_j(s)\,ds
\end{equation*}
$(r \ge (j-1)p)$.
Consequently,
\begin{align*}
|M_j(\lambda) - M_{j-1}(\lambda)| &\le |M_{j-1}(\lambda)| \int_{(j-1)p}^{jp}
|\Psi_{j-1}(s-p)^{-1}|\, \frac{|k|p}{s(s-p)}\,|\Psi_j(s)|\,ds \\
&\le C^3 |k|p\int_{(j-1)p}^{jp} \frac{ds}{s(s-p)}
\end{align*}
with a constant $C$ which is independent of $\lambda$ due to the uniform bound
(\ref{eq:psiphi}).
This also implies such an estimate for $|D_j(\lambda) - D_{j-1}(\lambda)|$ and,
since the $D_j$ are Lipschitz continuous on $\Lambda$, for
$|\mu_j(\lambda) - \mu_{j-1}(\lambda)|$ $(\lambda \in \Lambda)$.
Hence
\begin{equation*}
|E_j(\lambda) - E_{j-1}(\lambda)| \le C' \int_{(j-1)p}^{jp} \frac{ds}{s(s-p)}
\end{equation*}
with some other uniform constant $C'$.
Now we can estimate
\begin{align*}
|\Psi_J(np)| &\le |E_n|\,|E_J^{-1}| \prod_{j=J+1}^n |E_j^{-1} E_{j-1}|
\le C'' \prod_{j=J+1}^n (1 + |E_j^{-1}|\,|E_j - E_{j-1}|) \\
&\le C'' \exp \bigg(\sum_{j=J+1}^n |E_j^{-1}|\,|E_j - E_{j-1}|\bigg)
\le C'' \exp \left(C''' \int_{Jp}^{np} \frac{ds}{s(s-p)}\right),
\end{align*}
with uniform constants $C'', C'''$; this is bounded as $n \rightarrow \infty$.
Hence $\Psi_J(r)$ is bounded at infinity,
since $\Psi_J(r)= \Psi_n(r) \Psi_J((n-1)p)$ and, by (\ref{eq:psiphi}),
\begin{equation*}
|\Psi_n(r)| \le \big( 1 + \frac{1}{n-1} \big)^{|k|} \le 2^{|k|}  
\end{equation*}
for $(n-1)p \le r \le np$ and $n > J$.

This concludes the proof of Theorem \ref{th:radper}, since every solution of
(\ref{eq:raddir}) is a linear combination of the columns of $\Psi_J$.
$\;\blacksquare$

\medskip
The above method of proof does not work at the points
$\lambda\in\frac{\pi}{p}{\mathbb Z} + \hat\eta$;
these points are potential candidates for embedded eigenvalues.
However, it seems to be a rather delicate question to decide whether such embedded
eigenvalues actually occur.

\medskip

If $q$ is not assumed to be rotationally symmetric, we can prove 
the following. 

\vspace{3pt}

\begin{thm}  \label{th:virial}
Let $q \in C^1(\mathbb R^d; \mathbb R)$, and suppose
that both $q$ and $(x \cdot \nabla)q$ are 
bounded functions. 
Then  
$\sigma_p (H_d) \subset [m_q, \, M_q]$, where
\begin{equation*}
m_q= \inf_x \{ q(x) + (x\cdot \nabla)q(x) \}, 
\quad 
M_q= \sup_x \{ q(x) + (x\cdot \nabla)q(x) \}.
\end{equation*}
\end{thm}

\vspace{5pt}
To prove Theorem \ref{th:virial}, we shall apply  the following simple abstract version
of the virial theorem.

\begin{lem}[\textbf{Balinsky and Evans\cite{BE}}, \cite{BE-1}]
Let $U(a)$, $a>0$, be a one-parameter family of unitary operators
on a Hilbert space $\mathcal H$, which converges strongly to 
the identity as $a \to 1$. Let $T$ be
a self-adjoint operator in $\mathcal H$ and
$T_a:= aU(a)TU(a)^{-1}$.
If $f$ belongs to
$\text{\rm Dom}(T) \cap \text{\rm Dom}(T_a)$
and is an eigenvector of $T$ corresponding to 
an eigenvalue $\lambda$,
then
\begin{equation*}
\lim_{a \to 1} \Big(
f_a, 
\Big[
\frac{T_a - T}{a-1}
\Big] f
\Big)_{\mathcal H} 
=
\lambda \, \Vert f \Vert_{\mathcal H}^2,
\end{equation*}
where $f_a = U(a)f$.
\end{lem}

\noindent
{\bf Proof of Theorem \ref{th:virial}.}
We only give the proof for $H_3$, because the proof for $H_2$ is exactly the same.

Let $\lambda \in \sigma_p(H_3)$, and let
$f$ be a corresponding eigenfunction with $\Vert f \Vert =1.$
In particular, 
$f \in \mbox{Dom}(H_3)={\mathsf H}^1({\mathbb R}^3; {\mathbb C}^4)$
and 
 $H_3 f= \lambda f$. With the dilation group
$\{U(a)\}_{a>0}$, defined by 
$U(a)g (x):=a^{3/2}g(ax)$, $g \in {\mathsf L}^2({\mathbb R}^3; {\mathbb C}^4)$, 
we introduce a family of self-adjoint operators
$\{ H(a)\}_{a>0}$ by 
$H(a):= a \, U(a) H U(a)^{-1}$. 
We then see that
\begin{equation*}
\big( U(a)f, \, H(a)f \big) 
=  \big( H(a)U(a)f, \, f \big)
=  \big( aU(a) Hf, \, f \big)
=   \lambda a \big( U(a) f, \, f \big),
\end{equation*}
which implies that
\begin{equation}  \label{eq:vrl1}
\big( U(a)f, \, H(a)f - Hf \big)  = \lambda (a -1) \big( U(a) f, \, f \big).
\end{equation}

On the other hand, we find that
\begin{equation*}
[H(a) g](x)= -i \alpha \cdot \nabla g(x)+ a \,q (a x) g(x), 
\quad \forall g \in {\mathsf H}^1({\mathbb R}^3; {\mathbb C}^4),
\end{equation*}
hence we have
\begin{equation} \label{eq:vrl2}
 [H(a)f](x) - [Hf](x)  = \{a \,q (a x)  - q(x)\}f(x).
\end{equation}
Combining (\ref{eq:vrl1}) with (\ref{eq:vrl2}) 
yields
\begin{equation} \label{eq:vrl3}
\big( U(a)f, \, \frac{a \, q(a\,\cdot) - q( \,\cdot \,)}{a-1}f  \,\big)  
= \lambda  \big( U(a) f, \,   f\, \big).
\end{equation}
Since $\mbox{s\,-}\lim_{a \to 1} U(a)=I$, (\ref{eq:vrl3}) implies, by
the Lebesgue dominated convergence theorem, that
\begin{equation}    \label{eq:vrl4}
\int_{{\mathbb R}^3}
 \big( 
q(x) +  (x \cdot \nabla q)(x) 
\big)
|f(x)|^2 \, dx = \lambda.
\end{equation}
The conclusion of the theorem follows from (\ref{eq:vrl4}).
$\;\blacksquare$


\section{Schnol's theorem}
In this section, we state and prove Schnol's theorem for $H_d$.
The idea of our proof is based on that
of \cite[p. 21, Theorem 2.9]{CFKS},
where Schnol's theorem for Schr\"odinger operators 
is established.
In the three-dimensional case,
our Schnol's theorem can be stated as follows:

\begin{thm}  \label{schnol3}
Let $q\in {\mathsf L}^2_{loc}({\mathbb R}^3, {\mathbb R})$,
and let $\lambda$ be a real number.
Suppose $f$ is a polynomially bounded measurable function 
on ${\mathbb R}^3$, not identically $0$, and satisfies
the equation 
\begin{equation}  \label{eq:eveq3}
(- i \alpha \cdot \nabla + q) f = \lambda f
\end{equation}
in the distribution sense. 
Then $\lambda \in \sigma(H_3)$ 
for any
self-adjoint realization $H_3$ such that
$\mbox{\rm Dom}(H_3) 
\supset {\mathsf H}^1({\mathbb R}^3; {\mathbb C}^4)\cap \mbox{\rm Dom}(q)$.
\end{thm}

\vspace{2pt}

\noindent
{\bf Proof.} It is sufficient to prove the assertion for $\lambda =0$,
because any $\lambda\not=0$ can be absorbed in $q$. 
The proof will be devided into two steps.

\noindent
{\it Step 1. The case of $f \in {\mathsf L}^2({\mathbb R}^3; {\mathbb C}^4)$.} 

Let $\varphi \in C_0^{\infty}({\mathbb R}^3, {\mathbb C})$. Then
it follows that
\begin{align}
(- i\alpha \cdot \nabla )(\varphi f)
&=   
(- i\alpha \cdot \nabla \varphi) f 
+ 
\varphi (- i\alpha \cdot \nabla ) f   
 \nonumber  \\
&=  
(- i\alpha \cdot \nabla \varphi) f - \varphi \, q f,
\label{eqn:sa3-3}
\end{align}
where we have used (\ref{eq:eveq3})  in the second equality.
Since $\varphi q\in {\mathsf L}^2({\mathbb R}^3; {\mathbb C})$
and $f$ is locally bounded, 
we see  that 
$\varphi \, q f\in {\mathsf L}^2({\mathbb R}^3; {\mathbb C}^4)$, 
hence by (\ref{eqn:sa3-3}) that 
$(- i\alpha \cdot \nabla )(\varphi f)\in {\mathsf L}^2({\mathbb R}^3; {\mathbb C}^4)$.
This implies that 
$(- i\alpha \cdot \nabla )^2(\varphi f) \in 
{\mathsf H}^{-1}({\mathbb R}^3; {\mathbb C}^4)$, the Sobolev
space of order $-1$.
On the other hand, 
$(- i\alpha \cdot \nabla )^2(\varphi f)= -\Delta(\varphi f)$. 
Hence we find that $\varphi f \in 
{\mathsf L}^2({\mathbb R}^3; {\mathbb C}^4)
\subset 
{\mathsf H}^{-1}({\mathbb R}^3; {\mathbb C}^4)$
and $ -\Delta(\varphi f) \in  {\mathsf H}^{-1}({\mathbb R}^3; {\mathbb C}^4)$.
We now apply the ellipticity argument, 
and we get  
\begin{equation}  \label{eqn:sa2}
\varphi f \in {\mathsf H}^{1}({\mathbb R}^3; {\mathbb C}^4)
 \ \  \mbox{ for } \ 
\forall  \varphi \in C_0^{\infty}({\mathbb R}^3; {\mathbb C}).
\end{equation}
By the ellipticity argument, we mean the following:
\lq\lq 
$\, u \in {\mathsf H}^{\ell}({\mathbb R}^3; {\mathbb C})$ and 
$\Delta u \in {\mathsf H}^{\ell}({\mathbb R}^3; {\mathbb C})$ 
for some $\ell \in \mathbb R
\Longrightarrow  u \in 
{\mathsf H}^{\ell+2}({\mathbb R}^3; {\mathbb C}).$\rq\rq

We now choose $\chi \in C_0^{\infty}({\mathbb R}^3; {\mathbb C})$ 
so that
$\chi(x)=1 \ (|x| \le 1)$ and $=0 \ (|x| \ge 2)$,
and we set
\begin{equation}  \label{eqn:sa2-1}
\chi_n(x) = \chi(x/n)  \quad (n =1, \, 2, \, 3, \, \cdots).
\end{equation}
It follows from (\ref{eqn:sa2}) that 
$\chi_n f \in  {\mathsf H}^{1}({\mathbb R}^3; {\mathbb C}^4)$. 
It is evident that $\chi_n f \in \text{Dom}(Q)$.
Hence $\chi_n f \in \text{Dom}(H_3)$ for
$n =1, \, 2, \, 3, \, \cdots$.
To construct a singular sequence, we define
\begin{equation}  \label{eqn:sa3-5}
f_n = \frac{1}{\, \Vert \chi_n f \Vert_{{\mathsf L}^{2}}}\,
\chi_n f.
\end{equation}
It is obvious that 
$f_n \in \text{Dom}(H_3)$
and $\Vert f_n \Vert_{ {\mathsf L}^{2} }=1$.
We now only have to show that $\Vert H f_n \Vert_{{\mathsf L}^{2}}\to 0$
as $n \to \infty$.
In fact, we see that
\begin{gather}  \label{eqn:sa3-6}
\begin{split}
Hf_n
&=
 \frac{1}{\, \Vert \chi_n f \Vert_{{\mathsf L}^{2}}  }\,
\big[
\{ (-i \alpha \cdot \nabla )\chi_n \}f + 
\chi_n ( -i \alpha \cdot \nabla + q)f \big] \\
&=
\frac{1}{\, \Vert \chi_n f \Vert_{{\mathsf L}^{2}}  }\,
\Big[
\frac{1}{\, n \,}
\Big\{ 
(-i \alpha \cdot \nabla \chi) \big(\frac{x}{n} \big) 
\Big\} f  
\Big],
\end{split}
\end{gather}
where we have used the hypothesis   that 
$(-i \alpha \cdot \nabla + q) f =0$.
Noting the fact that 
$\lim_{n\to \infty}\Vert \chi_n f \Vert_{{\mathsf L}^{2}}
=\Vert f \Vert_{{\mathsf L}^{2}}\not= 0$,
we can deduce 
from (\ref{eqn:sa3-6})
that $\Vert H f_n \Vert_{{\mathsf L}^{2}}\to 0$.
Hence we can conclude that $0 \in \sigma(H_3)$.

\vspace{3pt}

\noindent
{\it Step 2. The case of $f \not\in {\mathsf L}^2({\mathbb R}^3; {\mathbb C}^4)$.}

We may assume that $f$ satisfies the estimate
\begin{equation}  \label{eqn:arara-1}
|f(x)| \le C (1+|x|)^N 
\end{equation}
for some $N \in \mathbb N$.
Let $f_n$ be defined in the same way as (\ref{eqn:sa3-5}). 
As was shown in Step 1, 
it follows that
$ f_n \in \mbox{\rm Dom}(H_3)$ and that (\ref{eqn:sa3-6}) is still valid.
Then we have
\begin{equation}    \label{eqn:sa3-8}
\Vert H f_n \Vert^2_{{\mathsf L}^2}
\le
\frac{1}{n^2 \Vert \chi_n f \Vert^2_{{\mathsf L}^2}} \,
 \Big[ \sup_{1\le |x| \le 2} |\nabla \chi (x)|  \Big]^2 
\int_{n \le |x| \le 2n} |f(x)|^2 \,dx.
\end{equation}

We now introduce a sequence
$(M(n))_{n \in {\mathbb N}}$ by
\begin{equation}   \label{eqn:sa3-9}
M(n):= \int_{|x| \le n} |f(x)|^2 \,dx
\quad (n\in {\mathbb N}),
\end{equation}
which is diverging and monotonically increasing.
It follows from (\ref{eqn:sa3-8}) and (\ref{eqn:sa3-9})
that
\begin{equation}  \label{eqn:arere-2}
\Vert H f_n \Vert^2_{{\mathsf L}^2}
\le
C \frac{M(2n) - M(n)}{n^2 M(n)},
\end{equation}
where $C$ is a positive constant, independent of $n$.

For the sake of contradiction, suppose that
\begin{equation}   \label{eqn:sa3-10}
\liminf_{n \to \infty}  \,\frac{M(2n) - M(n)}{n^2 M(n)} >0.
\end{equation}
Then there would be a large integer $\nu_0$ and
a positive constant $\alpha_0$ such that
\begin{equation}   \label{eqn:sa3-11}
\frac{M(2n) - M(n)}{n^2 M(n)} \ge \alpha_0
\quad \mbox{for } \forall n \ge \nu_0.
\end{equation}
This implies that
\begin{equation}   \label{eqn:sa3-12}
M(2n) \ge (1 + \alpha_0 n^2) M(n)
\quad \mbox{for } \forall n \ge \nu_0.
\end{equation}
By repeated use of (\ref{eqn:sa3-12}), we obtain
\begin{equation}   \label{eqn:sa3-13}
M(2^n \nu_0) \ge 
\Big\{
\prod_{j=0}^{n-1} (1 + \alpha_0 \nu_0^2 2^{2j}) 
\Big\}
M(\nu_0)
\quad (n\in{\mathbb N}).
\end{equation}
We now write $n=2\ell$. It follows from
(\ref{eqn:sa3-13}) that
\begin{gather}  \label{eqn:sa3-14}
\begin{split}
M(4^{\ell} \nu_0) 
&
\ge 
\Big\{
\prod_{j=0}^{2\ell-1} (1 + \alpha_0 \nu_0^2 4^j) 
\Big\}
M(\nu_0)  \\
&
\ge
\Big\{
\prod_{j=\ell}^{2\ell-1} (1+ \alpha_0 \nu_0^2 4^j) 
\Big\}
M(\nu_0)  \\
&
\ge
\alpha_0^{\ell} \nu_0^{2\ell} 4^{\ell^2}
M(\nu_0).
\end{split}
\end{gather}
On the other hand, it follows from (\ref{eqn:arara-1}) and
(\ref{eqn:sa3-9}) that
\begin{gather}   \label{eqn:sa3-15}
\begin{split}
M(4^{\ell} \nu_0) 
&
\le 
C \int_{|x| \le 4^{\ell} \nu_0} (1+|x|)^{2N} \,dx   \\
&
=
C^{\prime}
 \int_0^{4^{\ell} \nu_0} (1+r)^{2N} r^2 \, dr   \\ &
\le
C^{\prime\prime}_{\nu_0} 4^{(2N+3)\ell}.
\end{split}
\end{gather}
It follows from (\ref{eqn:sa3-14}) and (\ref{eqn:sa3-15})
that
\begin{equation}  \label{eqn:sa3-16}
C^{\prime\prime}_{\nu_0} 4^{(2N+3)\ell}  \ge 
\alpha_0^{\ell} \nu_0^{2\ell} 4^{\ell^2} M(\nu_0)
\quad (\ell \in{\mathbb N}).
\end{equation}
Taking the logarithm of both sides of (\ref{eqn:sa3-16}), 
one gets
\begin{gather}  \label{eqn:sa3-17}
\begin{split}
{}& \log C^{\prime\prime}_{\nu_0}
 + {(2N+3)\ell} \,\log 4  \ge  \\
{}& \qquad \qquad
\ell \,\log (\alpha_0 \nu_0^{2})
+  {\ell^2} \,\log4 + \log M(\nu_0) 
\qquad (\ell \in {\mathbb N}).
\end{split}
\end{gather}
Since the right hand side of (\ref{eqn:sa3-17})
grows faster than the left hand side of (\ref{eqn:sa3-17})
as $\ell$ goes to infinity,
the inequality (\ref{eqn:sa3-17})
 is a contradiction. Hence we can deduce
that
\begin{equation}   \label{eqn:sa3-18}
\liminf_{n \to \infty}  \,\frac{M(2n) - M(n)}{n^2 M(n)} =0,
\end{equation}
which yields that 
there is a subsequence $\{ M(n_k) \}$
 such that
\begin{equation}   \label{eqn:sa3-19}
\lim_{k \to \infty}  \,\frac{M(2n_k) - M(n_k)}{n_k^2 M(n_k)} =0.
\end{equation}
This fact, together with (\ref{eqn:arere-2}), implies that
$\Vert H f_{n_k} \Vert_{{\mathsf L}^2} \to 0$ 
as $k \to \infty$.
Thus we can conclude that $0 \in \sigma(H_3)$.
$\blacksquare$

\vspace{10pt}

In the two dimensional case,
Schnol's theorem is as follows:

 \begin{thm}  \label{schnol2}
 Let $q\in {\mathsf L}^2_{loc}({\mathbb R}^2, {\mathbb R})$  
and $\lambda$ be a real number.
Suppose $\psi$ is a polynomially bounded measurable function 
 on ${\mathbb R}^2$, not identically $0$, and satisfies
 the equation
 \begin{equation*}
 (- i \sigma \cdot \nabla + q) \psi = \lambda\psi
 \end{equation*}
 in the distribution sense. 
 Then $\lambda \in \sigma(H_2)$ 
 for any
 self-adjoint realization $H_2$ such that
 $\mbox{\rm Dom}(H_2) 
 \supset {\mathsf H}^1({\mathbb R}^2; {\mathbb C}^2)\cap \mbox{\rm Dom}(q)$.
 \end{thm}

The proof of Theorem \ref{schnol2} is similar to
 that of Thorem \ref{schnol3}, and is omitted.

When applying either of Theorems \ref{schnol3} and  \ref{schnol2},
one needs to construct a polynomially bounded eigensolution for
a given energy  of the Dirac operator $H_d$ 
with $q$ being locally square integrable. 
However, it is not easy to construct such an eigensolution
unless $q$ decays at infinity in an appropriate sense.
If $q$ decays rapidly, it is well-known that 
one can construct bounded eigensolutions (generalized eigenfunctions)
of $H_d$ by exploiting the limiting absorption principle. 
In Example \ref{schnlapp} below, 
we shall construct a bounded eigensolution (cf. (\ref{eq:sa-4-}))
 for a given energy of $H_d$ with potential $q$  of a specific form. 
We would like to stress that 
we do not require any decay assumption of $q$ at infinity.

\begin{expl} \label{schnlapp}
Let $\eta$ be a real-valued continuous function on $\mathbb R$ and define
$q(x):= \eta(x \cdot k)$ on ${\mathbb R}^d$, $d\in \!\{2, \, 3\}$,
where $k \in {\mathbb R}^d$ is a unit vector.
One can show, by the standard technique 
{\rm (}{\it cf.} {\rm \cite[{\it p.257, Corollary}\,]{RS-I}},
{\rm \cite[{\it p.113,  Theorem 4.3}\,]{Thaller}\rm )} ,
that $H_d$ is essentially self-adjoint on 
$H^1({\mathbb R}^d; {\mathbb C}^{2^{d-1}})\cap \mbox{\rm Dom}(q)$.
Let $H_d$ be the unique self-adjoint realization.
Then $\sigma(H_d)= \mathbb R$.

For $d=3$, this fact is proved in the following manner\,{\rm :} 
Put
\begin{equation*}
\xi(t) = \int_0^t \eta (\tau) \, d \tau.
\end{equation*}
As the eigenvalues of the matrix $\alpha\cdot k$ are $\pm 1$
 {\rm(}each with geometric multiplicity {\rm 2)}
we can choose a spinor  $\phi_0 \in {\mathbb C}^4$ 
so that $|\phi_0|=1$, $(\alpha \cdot k) \phi_0= \phi_0$.
For a given $\lambda \in \mathbb R$, define
\begin{equation}   \label{eq:sa-4-}
f(x)= e^{-i(\alpha\cdot k) \, \xi (x \cdot k)} 
 e^{i\lambda x\cdot k} \phi_0.
\end{equation}
Then $f$ is in $ C^1({\mathbb R}^3; {\mathbb C}^4)$, and
 satisfies the equation (\ref{eq:eveq3}).
In fact, one can see that
\begin{equation*}
-i \alpha_j \frac{\partial f}{\partial x_j}
=
-i \alpha_j  k_j  \, (-i \alpha \cdot k) \, q(x)   f  + \lambda \alpha_j k_j f
\qquad
(\, j\in \{1, \, 2, \, 3\}),
\end{equation*}
hence that 
\begin{align}
(-i \alpha\cdot \nabla )f
&=
 (-i \alpha \cdot k)^2 \, q(x)  f  + \lambda (\alpha \cdot k) f  \nonumber\\
&=
-q(x)  f  + \lambda  \, e^{-i(\alpha\cdot k) \, \xi (x \cdot k)} 
 e^{i\lambda x\cdot k} (\alpha \cdot k) \phi_0     \label{eq:sa-4}\\
&=
 -q(x)  f  + \lambda f,  \nonumber
\end{align}
where we have used the fact that $(\alpha \cdot k) \phi_0= \phi_0$
and the  facts that $(\alpha \cdot k)^2=I_4$ and
 that  $\alpha\cdot k$ commutes with the exponential 
$e^{-i(\alpha\cdot k) \, \xi (x \cdot k)}$.
It is obvious that 
 $| f(x) |_{{\mathbb C}^4} =1$ 
for all $x \in {\mathbb R}^3$.
Hence, it follows from 
Theorem \ref{schnol3} that
$\lambda \in \sigma(H_3)$.
Since $\lambda$ is an arbitrary real number, one can conclude that 
$\sigma(H_3)=\mathbb R$.

For $d=2$, the proof is similar to that for $d=3$ and is omitted.
\end{expl}
%


\section{Spectra of \boldmath $H_d$}

In this section, we shall prove that $\sigma(H_3)= \mathbb R$ under 
minimal assumptions on the potential q.
As mentioned in the introduction, the one-dimensional Dirac operator $H_1$
in (\ref{eq:do1d}) has (purely absolutely continuous) spectrum 
$\sigma(H_1)= \mathbb R$
for all potentials $q \in {\mathsf L}^1_{loc}({\mathbb R}; {\mathbb R})$.
In view of this fact,
the question naturally arises whether
the spectrum of $H_d$, $d \in \{2, \, 3\}$, also covers the whole 
real line
for all potentials $q \in {\mathsf L}^2_{loc}({\mathbb R}^d; {\mathbb R})$?
(Here we assume local square-integrability of the potential to ensure
that the Dirac operator will be well-defined on $C_0^\infty(\mathbb{R}^d; \mathbb{R})$.)

While attempting to give an answer to this question in greatest possible
generality,
we shall, however, need to impose some hypotheses on 
the potential $q$.
The reason for this restriction is technical.
Compared to the number of tools available to study the spectrum of
the one-dimensional Dirac operator (an ordinary differential operator),
the techniques for showing that a real number belongs to $\sigma(H_d)$
are relatively limited.
The main tool available for addressing the general question above is
Weyl's criterion in some form, {\it i.e.\/} the construction of a Weyl singular
sequence, as we have done in the previous section; and we shall use
this method again here.
However, we emphasize that the conditions we shall impose are fairly mild in that
they restrict the potential $q$ only on some sequence of balls of increasing radius,
which can be arbitrarily positioned and far apart.
In the remaining space, there is no constraint at all beyond the general assumption
of local square-integrability.
More specifically, in Theorem \ref{th:singseq3}, we need 
the potential $q$ to be sufficiently close, in
an ${\mathsf L}^2$ sense, to a function which varies only in one direction, 
and hence is constant on the planes perpendicular to this direction 
in each ball.
Theorem \ref{th:distorted}  is a generalization of Theorem \ref{th:singseq3},
where the planes of constancy are replaced with more general curved manifolds.
In  Theorem \ref{th:oscillation}, we need 
the mean oscillations of $q$ on the sequence of balls to go to zero.
This condition will be satisfied whenever the potential is close to
constant on wide stretches, even if these lie in e.g. a narrow
sector or cone. This indicates that a spectral gap could, if at all, only occur
in the case of a potential which changes in a complicated
multidimensional way essentially everywhere; an arrangement which would seem difficult to realise in practice.

Note that
we don't need any growth 
or decay property of the potential $q$ at infinity.

\begin{thm}  \label{th:singseq3}
Let $q \in {\mathsf L}^2_{loc}({\mathbb R}^3; {\mathbb R})$. 
Suppose that there is a sequence $(k_n)_{n\in{\mathbb N}}$ of
unit vectors in ${\mathbb R}^3$,
a sequence $(B_{r_n}(a_n))_{n\in{\mathbb N}}$ of 
disjoint
balls with centre
$a_n \in {\mathbb R}^3$ and radius $r_n \to \infty$ $(n\rightarrow\infty)$,
and a sequence of square-integrable functions
$\eta_n: (-r_n,r_n)\rightarrow{\mathbb R}$ $(n\in{\mathbb N})$ such  that
\begin{equation*} 
r_n^{-3}
\int_{B_{r_n}(a_n)} 
 \big|q(x) - \eta_n \big( (x - a_n)\cdot k_n) \big) \big|^2
\, dx \to 0 
\end{equation*}
as $n \to \infty$. Then $\sigma(H_3) = {\mathbb R}\;$   
for
any self-adjoint extension $H_3$ of
\begin{equation*}  
(-i\alpha \cdot \nabla + q) 
\big|_{C_0^\infty({\mathbb R}^3;  {\mathbb C}^4)}.
\end{equation*}
\end{thm}

\vspace{6pt}

\begin{rem}
The two dimensional analogue of Theorem \ref{th:singseq3} holds true.
\end{rem}

\noindent
{\bf Proof of Theorem 4.1.}
Let $\lambda\in {\mathbb R}$; we shall show that $\lambda$ belongs to the spectrum
of $H_3$ by constructing a Weyl singular sequence.

In a similar way to Example \ref{schnlapp}, we can choose a sequence of
spinors $(\phi_n)_{n\in{\mathbb N}}$ in ${\mathbb C}^4$ such that
\begin{equation} \label{eq:phi-n}
|\phi_n|=1, \quad (\alpha \cdot k_n) \phi_n = \phi_n.
\end{equation}
Since $C^\infty(-r_n, r_n)$ is dense in ${\mathsf L}^2(-r_n, r_n)$,
there are functions $\tilde \eta_n\in C^\infty(-r_n, r_n)$ such that
\begin{equation}
\frac{1}{2 r_n}\int_{-r_n}^{r_n} |\eta_n(\tau) - \tilde \eta_n(\tau)|^2  d\tau \to 0
\end{equation}
 as
 $n\to\infty$.
Let $\xi_n(t) := \int_0^t \tilde \eta_n (\tau) d\tau$ \ 
$(t\in(-r_n, r_n); n\in{\mathbb N})$, and
define
\begin{equation}  \label{eq:Fn}
F_n(x) := e^{-i(\alpha \cdot k_n)
\xi_n((x-a_n)\cdot k_n)}
e^{i \lambda x \cdot k_n} \phi_n
\, : \, B_{r_n}(a_n) \to {\mathbb C}^4.
\end{equation}
Since $e^{-i(\alpha \cdot k_n)\xi_n((x-a_n)\cdot k_n)}$ is 
a unitary matrix, it follows
from (\ref{eq:phi-n}) that 
$|F_n(x)|_{{\mathbb C}^4}=1$ for all $n\in{\mathbb N}$ and all $x \in B_{r_n}(a_n)$.
Furthermore, we see that $F_n \in C^\infty(B_{r_n}(a_n))^4$, and
we make the same computation as in (\ref{eq:sa-4}) to get
\begin{gather} \label{eq:DOFn}
\begin{split}
(-i \alpha \cdot \nabla) F_n(x)
=
\big\{ \! 
- \tilde \eta_n\big((x - a_n) \cdot k_n\big) + \lambda  \,
\big\}\,F_n(x).
\end{split}
\end{gather}
Here we have used (\ref{eq:phi-n})
and the facts that $(\alpha\cdot k_n)^2 = I_4$ and that $\alpha\cdot k_n$
commutes with the exponentials in (\ref{eq:Fn}).

We now choose $\chi \in C_0^{\infty}({\mathbb R}^3; {\mathbb R})$ 
so that
$\mbox{supp} (\chi) \subset B_1(0)$ and
$\Vert \chi \Vert_{\mathsf L^2}=1$, and 
define
\begin{equation}  \label{eq:cut}
\chi_n(x):= r_n^{-3/2} \, 
\chi\big( r_n^{-1}(x - a_n)\big).
\end{equation}
We shall show that the sequence $(f_n)_{n\in{\mathbb N}}$ defined
by $f_n:=\chi_n F_n$ $(n\in{\mathbb N})$ is a Weyl singular sequence for
$H_3-\lambda$.
First, we note that $\Vert f_n \Vert_{{\mathsf L}^2}=1$
and $f_n \in C_0^\infty({\mathbb R}^3;{\mathbb C}^4)$.
Next, we see that
\begin{gather}  \label{eq:DOfn}
\begin{split}
(-i \alpha \cdot \nabla) f_n(x) 
&=
\big\{(-i \alpha \cdot \nabla)\chi_n(x) \big\}\,F_n(x) 
+
\chi_n(x)(-i \alpha \cdot \nabla) F_n(x)   \\
&=  r_n^{-5/2} \, 
\big[
(-i \alpha \cdot \nabla)\chi
\big]
\big( r_n^{-1}(x - a_n)\big) F_n(x) \\
&\qquad 
-  \tilde \eta_n \big((x - a_n) \cdot k_n \big)\,f_n(x)
+ \lambda f_n(x),
\end{split}
\end{gather}
where we have used (\ref{eq:DOFn}) and (\ref{eq:cut}).
Finally, it follows from (\ref{eq:DOfn})
that
\begin{gather}  \label{eq:DOfn1}
\begin{split}
(H_3-\lambda)f_n(x)
&=  
r_n^{-5/2} \, \big[(-i \alpha \cdot \nabla)\chi\big]
\big( r_n^{-1}(x - a_n)\big) \, F_n(x)  \\
{}&\;\; +
\big\{
q(x)- \tilde \eta_n\big((x - a_n) \cdot k_n\big) 
\big\}
 \, r_n^{-3/2} \, 
\chi\big(r_n^{-1}(x - a_n) \big) \, F_n(x) ,
\end{split}
\end{gather}
which --- adding and subtracting 
$\eta_n\big((x - a_n)\cdot k_n\big)$ --- implies that
\begin{gather} \label{eq:DOfn2}
\begin{split}
\Vert (H_3-\lambda)f_n \Vert
&\le 
r_n^{-1} \, \left(
\int_{|x|\le 1}
 \big| (\alpha \cdot \nabla)\chi(x) \big|^2
\, dx
\right)^{1/2}  \\
\noalign{\vskip 4pt}
&\;\;+
\Vert \chi \Vert_\infty \,
\left(
 r_n^{-3} 
\int_{B_{r_n}(a_n)}
\big| q(x)- \eta_n\big( (x - a_n) \cdot k_n \big) \big|^2
 \,  dx  \right)^{1/2}   \\
&\;\;+
\Vert \chi \Vert_\infty \,
\left(
\frac{\pi}{r_n} 
\int_{-r_n}^{r_n}
\big|\eta_n(\tau)-\tilde \eta_n(\tau) \big|^2 d\tau \right)^{1/2}   \\
\noalign{\vskip 4pt}
&\quad \to 0
\ \  \mbox{ as } \  n \to \infty.
\end{split}
\end{gather}
Hence we can conclude that $\lambda \in \sigma(H_3)$. 
$\quad\blacksquare$

\vspace{8pt}

In the following theorem,
 we shall show that the result of Theorem \ref{th:singseq3}
extends to potentials which are
close to  constants on a local foliation of curved surfaces, which could
be fattened to sets of positive measure, 
provided  that 
their curvature becomes
asymptotically small. 
We assume that the potential is
approximated by  $C^\infty$ smooth functions,
which in the light of the proof of the preceding theorem is no serious
restriction of generality.

\vspace{10pt}

\begin{thm}  \label{th:distorted}
Let $q \in {\mathsf L}^2_{loc}({\mathbb R}^3; {\mathbb R})$. 
Suppose that there is a sequence $(k_n)_{n\in{\mathbb N}}$ of
unit vectors in ${\mathbb R}^3$,
a sequence $(B_{r_n}(a_n))_{n\in{\mathbb N}}$ of 
disjoint
balls with centre
$a_n \in {\mathbb R}^3$ and radius $r_n \to \infty$ $(n\to\infty)$,
and sequences of functions
$\varphi_n \in C^\infty(B_{r_n}(a_n); {\mathbb R})$
and
$\eta_n \in C^\infty({\mathbb R};{\mathbb R})$ $(n\in{\mathbb N})$
such  that
\begin{equation} \label{eq:phin}
r_n^{-3}
\int_{B_{r_n}(a_n)} 
 \big|q(x) - \eta_n \big(x\cdot k_n + \varphi_n(x) \big) \big|^2
\, dx \to 0  
\quad (n \to \infty)
\end{equation}
and
\begin{equation} \label{eq:phin1}
r_n^{-3}
\int_{B_{r_n}(a_n)} 
\big|\nabla \phi_n(x) \big|^2 \,
 \big|\eta_n \big(x\cdot k_n + \varphi_n(x) \big) \big|^2
\, dx \to 0  
\quad (n \to \infty).
\end{equation}
Then $\sigma(H_3) = {\mathbb R}$   
for any self-adjoint extension $H_3$ of
\begin{equation}  \label{eq:domain2}
(-i\alpha \cdot \nabla + q) 
\big|_{C_0^\infty({\mathbb R}^3;  {\mathbb C}^4)}
\end{equation}
\end{thm}

\vspace{6pt}

\begin{rem}
The two dimensional analogue of Theorem \ref{th:distorted} holds true.
\end{rem}

\noindent
{\bf Proof of Theorem 4.2.}
We follow the general line of the proof of 
Theorem \ref{th:singseq3}. Let $\lambda \in {\mathbb R}$ be arbitrary, and
$\phi_n$ as in (\ref{eq:phi-n}).
Define $\xi_n(t) := \int_0^t \eta_n (\tau) \, d\tau$ 
$(t\in{\mathbb R})$ and
\begin{equation}  \label{eq:Fnphi}
F_n(x) := 
e^{-i(\alpha \cdot k_n)\xi_n(x\cdot k_n + \varphi_n(x))}
e^{i\lambda\, x \cdot k_n} \phi_n
\qquad
(x\in B_{r_n}(a_n));
\end{equation}
Then $F_n \in C^\infty(B_{r_n}(a_n), {\mathbb C}^4)$ and
$|F_n(x)|_{{\mathbb C}^4}=1$ $(n\in{\mathbb N}, x \in B_{r_n}(a_n))$.
Moreover, abbreviating $q_n(x) := \eta_n(x\cdot k_n + \varphi_n(x))$
we get
\begin{gather}
\begin{split}
-i\alpha\cdot\nabla F_n(x)
&=
-i\alpha\cdot\big(
-i(\alpha\cdot k_n)\,
q_n(x)\,
 \{k_n + \nabla\varphi_n(x) \} + i \lambda\,k_n
\big)\,F_n(x) \\
&=
-q_n(x)\,F_n(x)
 -\big(  \alpha\cdot\nabla\varphi_n(x)   \big)
\,(\alpha\cdot k_n)\,q_n(x)\,F_n(x)
 + \lambda\,F_n(x)
\end{split}
\end{gather}
$(x\in B_{r_n}(a_n))$.
To construct a singular sequence, let $\chi_n$ be as in (\ref{eq:cut}), and 
define $f_n:= \chi_n F_n$.
Then $\Vert f_n \Vert_{{\mathsf L}^2}=1$
and $f_n \in C_0^\infty({\mathbb R}^3; {\mathbb C}^4)$.
Furthermore, we have
\begin{gather}  \label{eq:DOfnphi}
\begin{split}
-i \alpha \cdot \nabla f_n(x) 
&=
\big( -i \alpha \cdot \nabla\chi_n(x)  \big) F_n(x) 
+
\chi_n(x) \big(  -i \alpha \cdot \nabla F_n(x)  \big)   \\
&=  r_n^{-5/2}\, \big[(-i \alpha \cdot \nabla)
\chi\big]
(r_n^{-1} \big(x - a_n) \big) \,F_n(x) \\
&\quad - \chi_n(x) \, q_n(x)\,F_n(x) \\
&\quad
-\chi_n(x)
\big(
\alpha \cdot \nabla \varphi_n(x)
\big) 
(\alpha \cdot k_n)
q_n(x)\,F_n(x) 
+ \lambda f_n(x),
\end{split}
\end{gather}
from which we obtain
\begin{gather}   \label{eq:singsq}
\begin{split}
(H_3-\lambda)f_n(x)
&=
r_n^{-5/2} \, 
\big[
(-i \alpha \cdot \nabla)\chi
\big]
\big( r_n^{-1}(x - a_n) \big)\, F_n(x) \\
&\;\; +
\chi_n(x) \big(q(x) - q_n(x) \big)\,F_n(x)  \\
&\;\;
-\chi_n(x)
\big( \alpha \cdot \nabla \varphi_n(x) \big) (\alpha \cdot k_n) \,
q_n(x) F_n(x).
\end{split}
\end{gather}
Hence
\begin{gather}   \label{eq:singsq1}
\begin{split}
\Vert  (H_3-\lambda)f_n \Vert_{{\mathsf L}^2}
&\le 
r_n^{-1} \, \Big(
\int_{|x|\le 1}
 \big| (\alpha \cdot \nabla)\chi(x) \big|^2
\, dx
\Big)^{1/2}  \\
\noalign{\vskip 4pt}
&+
\Vert \chi \Vert_{{\mathsf L}^{\infty}} \,
\Big(
 r_n^{-3} 
\int_{B_{r_n}(a_n)}
\big| q(x)- \eta_n\big( (x - a_n) \cdot k_n \big) \big|^2
 \,  dx  \Big)^{1/2}   \\
\noalign{\vskip 4pt}
&+
\Vert \chi \Vert_{{\mathsf L}^{\infty}} \,
\Big(
 r_n^{-3} 
\int_{B_{r_n}(a_n)}
\big| \nabla \varphi_n (x) \big|^2 \, 
\big| \eta_n\big( (x - a_n) \cdot k_n \big) \big|^2
 \,  dx  \Big)^{1/2}   \\
\noalign{\vskip 4pt}
&\quad \to 0
\ \  \mbox{ as } \  n \to \infty, 
\end{split}
\end{gather}
by (\ref{eq:phin}) and  (\ref{eq:phin1}).
Here we twice used the fact that 
$|({\alpha}\cdot v)u|=|v||u|$ for 
any $v \in {\mathbb R}^3$ and $u \in {\mathbb C}^4$. 
Thus we can conclude that $\lambda \in \sigma(H_3)$. 
$\quad\blacksquare$

\vspace{10pt}

The following theorem can be obtained as a special case of Theorem \ref{th:singseq3}, when
the functions $\eta_n$ are taken to be constants with value $q_n$ defined in
(\ref{eq:oscill2}).
However, it has a quick and simple separate proof which we include below.
\vspace{10pt}
\begin{thm}  \label{th:oscillation}
Let $q \in {\mathsf L}^2_{loc}({\mathbb R}^3; {\mathbb R})$. 
Suppose that there is 
a sequence $(B_{r_n}(a_n))_{n\in{\mathbb N}}$ of disjoint
balls with centre
$a_n \in {\mathbb R}^3$ and radius $r_n \to \infty$ $(n\to\infty)$
such  that
\begin{equation} \label{eq:oscill1}
r_n^{-3}
\int_{B_{r_n}(a_n)} 
 \big|q(x) - q_n \big|^2
\, dx \to 0  
\quad (n \to \infty),
\end{equation}
where $q_n$ is the mean value of $q$ over the ball $B_{r_n}(a_n)$:
\begin{equation} \label{eq:oscill2}
q_n=
 \frac{3}{\,4{\pi} r_n^{3}\,}\int_{B_{r_n}(a_n)} q(x) \, dx.
\end{equation}
Then $\sigma(H_3) = {\mathbb R}$   
for any self-adjoint extension $H_3$ of
\begin{equation*}  \label{eq:domain3}
(-i\alpha \cdot \nabla + q) 
\big|_{C_0^\infty({\mathbb R}^3;  {\mathbb C}^4)}
\end{equation*}
\end{thm}

\vspace{2pt}

\noindent
{\bf Proof.}  Let $\lambda \in {\mathbb R}$ be arbitrary.
As in Theorems \ref{th:singseq3} and \ref{th:distorted},
we shall show that $\lambda$ belongs to the spectrum
of $H_3$ by constructing a Weyl singular sequence.

We first choose a sequence $(k_n)_{n\in{\mathbb N}}$ in ${\mathbb R}^3$ so that 
$|k_n|=| q_n - \lambda|$ for each $n$. 
Note that $k_n=0$ if $q_n=\lambda$.
Since the eigenvalues of the matrix $\alpha\cdot k_n$ are $\pm |k_n|$,
we can choose a sequence of spinors  
$(\phi_n)_{n\in{\mathbb N}}$ in ${\mathbb C}^4$ so that
\begin{equation} \label{eq:phi-n-3}
|\phi_n|=1, \quad (\alpha \cdot k_n) \phi_n = 
\begin{cases}
-|k_n|\phi_n  &  \text{if } q_n - \lambda >0 \\
\;\;\, |k_n|\phi_n & \text{if } q_n - \lambda < 0 \\
\;\;\;\; 0 &  \text{if } q_n - \lambda =0,
\end{cases}
\end{equation}
which readily implies that 
\begin{equation}   \label{eq:phi-n-4}
(\alpha \cdot k_n + q_n - \lambda) \phi_n =0.
\end{equation}
With $\chi_n$ introduced in (\ref{eq:cut}), 
we shall show that the sequence $(f_n)_{n\in{\mathbb N}}$ defined
by $f_n:=\chi_n  \, e^{ix\cdot k_n}\phi_n$ is a Weyl singular sequence for
$H_3-\lambda$.
To this end, we first note that $\Vert f_n \Vert_{{\mathsf L}^2}=1$
and $f_n \in C_0^\infty({\mathbb R}^3;{\mathbb C}^4)$.
We next see that
\begin{gather}  \label{eq:DOf-n-3}
\begin{split}
(-i \alpha \cdot \nabla) f_n(x) 
&=  r_n^{-5/2} \, 
\big[
(-i \alpha \cdot \nabla)\chi
\big]
\big( r_n^{-1}(x - a_n)\big)  \ e^{ix\cdot k_n}\phi_n   \\
\noalign{\vskip 4pt}
&\;\;\;
+  \chi_n(x)  \, e^{ix\cdot k_n}(\alpha \cdot k_n)\phi_n
\end{split}
\end{gather}
Combining (\ref{eq:phi-n-4}) and (\ref{eq:DOf-n-3}),
we get
\begin{gather}  \label{eq:DOf-n-4}
\begin{split}
\big( \! -i \alpha \cdot \nabla + q(x) - \lambda \big) f_n(x) 
&=
r_n^{-5/2} \, 
\big[
(-i \alpha \cdot \nabla)\chi
\big]
\big( r_n^{-1}(x - a_n)\big)  \ e^{ix\cdot k_n}\phi_n   \\
\noalign{\vskip 4pt}
&\;\;\;
+  \chi_n(x)  \, e^{ix\cdot k_n}\big(q(x) - q_n \big)\phi_n,
\end{split}
\end{gather}
which implies
\begin{gather}   \label{eq:DOf-n-5}
\begin{split}
\Vert  (H_3-\lambda)f_n \Vert_{{\mathsf L}^2}
&\le 
r_n^{-1} \, \Big(
\int_{|x|\le 1}
 \big| (\alpha \cdot \nabla)\chi(x) \big|^2
\, dx
\Big)^{1/2}  \\
\noalign{\vskip 4pt}
&+
\Vert \chi \Vert_{{\mathsf L}^{\infty}} \,
\Big(
 r_n^{-3} 
\int_{B_{r_n}(a_n)}
\big| q(x)- q_n \big|^2
 \,  dx  \Big)^{1/2}  
\to 0
\end{split}
\end{gather}
as $n \to \infty$,
by assumption (\ref{eq:oscill1}).$\quad\blacksquare$

\vspace{6pt}

\begin{rem}
In (\ref{eq:oscill1}), the mean oscillation is taken in ${\mathsf L}^2$ sense.
The mean oscillation in the usual sense  is, however, taken in ${\mathsf L}^1$
sense; see {\it e.g.} \cite{Stein}.  One can see that (\ref{eq:oscill1}) implies the 
mean oscillation in the usual sense tends to zero as follows:
\begin{align*}
&\frac{1}{\,|B_{r_n}(a_n)| \,}
\int_{B_{r_n}(a_n)} 
 \big|q(x) - q_n \big|\, dx \\
\noalign{\vskip 4pt}
&\le
\frac{1}{\,|B_{r_n}(a_n)| \,}
\Big\{ \int_{B_{r_n}(a_n)} 
 \big|q(x) - q_n \big|^2 \, dx\Big\}^{1/2}
   \Big\{ \int_{B_{r_n}(a_n)} \, dx\Big\}^{1/2}  \\
\noalign{\vskip 4pt}
&=
\Big\{ 
\frac{3}{\, 4\pi \,} r_n^{-3}   
\int_{B_{r_n}(a_n)}   
 \big|q(x) - q_n \big|^2 \, dx\Big\}^{1/2}
\, dx \to 0 
\quad (n \to \infty),
\end{align*}
where $|\cdot|$ denotes the Lebesgue measure: $|B_{r_n}(a_n)|= 4\pi r_n^3/3$.
\end{rem}

\vspace{1pt}

\begin{rem}
The two dimensional analogue of Theorem \ref{th:oscillation} holds true.
\end{rem}

\vspace{1pt}

\noindent
{\bf Acknowledgements}

\noindent
The authors would like to thank the referee for the valuable comments
which led to an improvement of the paper.

\end{document}